\input amstex
\documentstyle{conm-p}
\NoBlackBoxes
\define\ra{\rightarrow}
\define\W{\omega}
\define\inc{\subseteq}

%  The following items provide publication information for the AMS-P logo
\issueinfo{00}% volume number
  {}% issue number
  {}% month
  {XXXX}% year

\topmatter
\title Hyman Bass and Ubiquity: Gorenstein Rings
\endtitle

\author Craig Huneke  \endauthor
\leftheadtext{CRAIG HUNEKE}

\address
Department of Mathematics, University of Kansas,
Lawrence, KS 66045 \endaddress
\email  huneke\@math.ukans.edu  \endemail

%  \thanks will become a 1st page footnote.
%  Use \endgraf to indicate a new paragraph; a blank line or \par will
%  be recognized as an error.
\thanks The author was supported in part by the NSF. I would like to
thank the many people who read preliminary versions and made
suggestions. In particular, I would like to thank Winfried Bruns
and the referee for valuable comments.
\endthanks

%  Math Subject Classifications
\subjclass Primary 13CDH ; Secondary 14B \endsubjclass

\dedicatory Dedicated to Hyman Bass.\enddedicatory

\keywords Gorenstein rings \endkeywords

\abstract
This paper is based on a talk given by the author in October, 1997
at a conference at Columbia University in celebration of Hyman
Bass's 65th birthday. The paper details some of
the history of Gorenstein rings and their uses.
\endabstract

\toc
\title Introduction
\page {1}
\endtitle
\title Plane Curves
\page{2}
\endtitle
\title Commutative Algebra circa 1960
\page{4}
\endtitle
\title Gorenstein Rings
\page{6}
\endtitle
\title
Examples and Low Codimension
\page{9}
\endtitle
\title
Injective Modules and Matlis Duality
\page{11}
\endtitle
\title
Ubiquity
\page{13}
\endtitle
\title
Homological Themes
\page{15}
\endtitle
\title
Inverse Powers and $0$-dimensional Gorenstein Rings
\page{16}
\endtitle
\title
Hilbert Functions
\page{18}
\endtitle
\title
Invariants and Gorenstein Rings
\page{20}
\endtitle
\title
Ubiquity and Module Theory
\page{21}
\endtitle
\title Bibliography
\page{23}
\endtitle
\endtoc
\endtopmatter
\document
\bigskip
\head  Introduction \endhead

In 1963, an article appeared in Mathematische Zeitschrift
with an interesting title, `On the ubiquity of Gorenstein rings', and
fascinating content. The
author was Hyman Bass. The article has become one of the most-read and
quoted math articles in the world.  A journal survey done in 1980 \cite{Ga}
showed that it ranked third among most-quoted papers from core
math journals. For every young student in commutative algebra, it is
at the top of a list of papers to read.
 The first thing I did as a graduate student when
I opened the journal to this paper was to reach for a dictionary. 

\definition{Definition: Ubiquity} The state or capacity of
being everywhere, especially at the same time; omnipresent.
\enddefinition

  This article will give some of the historical
background of Gorenstein rings, explain a little of why they are
ubiquitous and useful, and give practical ways of computing them.
The property of a ring being Gorenstein is fundamentally a
statement of symmetry. 
The `Gorenstein' of Gorenstein rings is Daniel Gorenstein, the
same who is famous for his role in the classification of finite
simple groups. A question occurring to everyone who studies Gorenstein
rings is `Why are they called Gorenstein rings?'
His name being attached to this concept goes
back to his thesis on plane curves, written under Oscar Zariski and published in
the Transactions of the American Mathematical Society in 1952 \cite{Go}.
As we shall see, they could perhaps more justifiably be called Bass rings, or  Grothendieck rings, 
or Rosenlicht rings, or Serre rings.   The usual definition
now used in most textbooks goes back to the work of Bass in the ubiquity
paper.  Going back even further, one could make on argument that
the origins of Gorenstein rings lie in the work of
W. Gr\" obner, and F.S. Macaulay. Indeed,
a 1934 paper of  Gr\" obner \cite{Gro}\footnote{I am very grateful to
the referee for pointing out this paper and its relevance to the
discussion of Gorenstein rings.} explicitly gives  the basic
duality of a $0$-dimensional Gorenstein ring and recognizes the
role of the socle: see Section 5 for a discussion of this duality.

\bigskip

\head 1. Plane Curves
\endhead
\bigskip
The origins of Gorenstein rings, at least as far as the history of
the definition of them,  go back to the classical study of
plane curves. Fix a field $k$
and let $f(X,Y)$ be a polynomial in the ring $k[X,Y]$. By a plane
curve we will mean the ring $R = k[X,Y]/(f)$. If $k$ is algebraically
closed, Hilbert's Nullstellensatz allows us to identify the maximal
ideals of this ring with the solutions of $f(X,Y) = 0$ which lie in $k$,
via the correspondence of solutions $(\alpha,\beta)\in k^2$ with 
maximal ideals $(X-\alpha, Y-\beta)$.
In his thesis, Gorenstein was interested in the properties of the
so-called adjoint curves to an irreducible plane curve $f = 0$. However,
his main theorem had to do with the integral closure of a plane curve.

\definition{Definition 1.1} Let $R$ be an integral domain with
fraction field $K$. The \it integral closure \rm of $R$ is
the set of all elements of $K$ satisfying a monic polynomial
with coefficients in $R$.
\enddefinition

The integral closure is a ring $T$, containing $R$, which is
itself integrally closed. An important measure of the
difference between $R$ and its integral closure $T$ is
the \it conductor ideal\rm, denoted $\frak C = \frak C(T,R)$,
$$
\frak C = \{r\in R|\, rT\subseteq R\}.
$$

It is the largest common ideal of both $R$ and $T$.

\example{Example 1.2} Let $R = k[t^3,t^7]\cong
k[X,Y]/(X^7-Y^3)$. The element 
$$t = \frac {t^7} {(t^3)^2},$$
is in the fraction field of $R$ and is clearly integral over $R$. It follows that the integral
closure of $R$ is $T = k[t]$. The dimension over $k$ of
$T/R$ can be computed by counting powers of  $t$ which are in $T$ but not in $R$,
since the powers of $t$ form a $k$-basis of both $R$ and $T$. The conductor is
the largest common ideal of both $T$ and $R$. It will contain some least
power of $t$, say $t^c$, and then must contain all higher powers of $t$
since it is an ideal in $T$. To describe $R$, $T$ and $\frak C$, it suffices
to give the exponents of the powers of $t$ inside each of them.
\bigskip

For $T$: $0, 1,\, 2, 3, 4, 5, \, 6, 7, 8, 9, 10, 11, 12, 13, ...$

For $R$: $0, \quad \quad 3,\quad \quad  6, 7, \quad 9, 10,\quad\,\,\, 12, 13,...$

For $\frak C$: $\quad \quad \quad \quad \quad \quad \quad \quad \quad \quad \quad \quad
\quad 12, 13$...
\endexample

A surprising phenomenon can be found by examining this chart; the number
of monomials in $T$ but not in $R$ is the same as the number in $R$ but
not in the conductor $\frak C$. In this example there are six such monomials.
More precisely, the vector space dimension of $T/R$ is equal to that
of $R/\frak C$.  Is this an accident? It is no accident, and this equality is
the main theorem of Gorenstein, locally on the curve. 
In fact, independently 
Ap\' ery \cite{A} in 1943, Samuel \cite{Sa} in 1951 and Gorenstein \cite{Go} in 1952 all
proved that given a prime $Q$ of a plane curve $R$ with integral closure
$T$ and conductor $\frak C$ 

$$
\text{dim}_k(T_Q/R_Q) = \text{dim}_k(R_Q/\frak C_Q)
$$
where $T_Q, R_Q$ and $\frak C_Q$ are the localizations at the prime $Q$.
(In other words we invert all elements outside of $Q$; $R_Q$ is then a
local ring with unique maximal ideal $Q_Q$. Recall that a \it local ring \rm
is a Noetherian commutative ring with
a unique maximal ideal $m$.) 
It turns out that \it all \rm $1$-dimensional local Gorenstein domains have this
fundamental equality; in fact it essentially characterizes such $1$-dimensional
domains. Plane curves are a simple  example of Gorenstein
rings; the simplest examples are polynomial rings or formal power series
rings. 

More general than plane curves are complete intersections.
Let $S$ a regular ring (e.g. the polynomial ring  $ k[X_1,...,X_n]$  or a
regular local ring) and let $R = S/I$. The codimension of $R$ (or $I$) is
defined by
$$
\text{codim}(R) = \text{dim}(S) - \text{dim}(R) = n - \text{dim}(R).
$$

$R$ is said to be a \it complete intersection \rm if $I$ can be generated by
$\text{codim}(R)$ elements.
Plane curves are a simple example since
they are defined by one equation. All complete intersections are Cohen-Macaulay,
and the defining ideal $I$ will be generated by a regular sequence.
An important example of a complete intersection is
$\bold Z[\Pi]$
where $\Pi$ is a finitely generated abelian group.

In 1952, Rosenlicht \cite{Ro} proved the equality
$$ \text{dim}_k(T_Q/R_Q) = \text{dim}_k(R_Q/\frak C_Q)$$
for localizations of complete
intersection curves, i.e. for one-dimensional domains $R$ of the form
$$ R = k[X_1,...,X_n]/(F_1,...,F_{n-1})
$$
with integral closure $T$ and conductor $\frak C$.

Furthermore, Rosenlicht proved that a local ring $R$ of an algebraic
curve had the property that $\text{dim}_k(T/R) = \text{dim}_k(R/\frak C)$
iff the module of regular differentials $\Omega$ is a free $R$-module.
It is interesting that both Gorenstein and Rosenlicht were students of
Zariski. As for plane curves, complete intersections are also special
examples of Gorenstein rings.
To explain the background and definitions, we need to review the
state of commutative algebra around 1960.

\bigskip

\head 2. Commutative Algebra circa 1960\endhead
\bigskip

The most basic commutative rings are the polynomials rings over a field
or the integers, and their quotient rings. Of particular importance are the rings
whose corresponding varieties are non-singular. To understand the
notion of non-singularity in terms of the ring structure,
the notion of regular local rings was developed.

\definition{Definition-Theorem 2.1} A \it regular local ring \rm $(S,n)$ is a Noetherian local ring which
satisfies one of the following equivalent properties:

\roster
\item The maximal ideal $n$ is generated by dim$(S)$ elements.
\item Every finitely generated $S$-module $M$ has a finite resolution
by finitely generated free $S$-modules, i.e. there is an exact sequence,
$$
0\ra F_k\overset {\alpha_k}\to\ra F_{k-1}\overset {\alpha_{k-1}}\to\ra ...\ra F_0\ra M
\ra 0,
$$
where the $F_i$ are finitely generated free $S$-modules.
\item The residue field $S/n$ has a finite free resolution as in (2).
\endroster
\enddefinition

These equivalences are due to Auslander-Buchsbaum \cite{AuBu1} and Serre \cite{Se2}, 
and marked a watershed in commutative algebra. One immediate use of this
theorem was to prove that
the localization of a regular local ring is still regular.
 Regular local rings correspond to
the localizations of affine varieties at non-singular points. The
existence of a finite free resolution for finitely generated modules
is of crucial importance. For example, in 1959 Auslander and Buchsbaum \cite{AuBu2}
used the existence of such resolutions to prove their celebrated result
that regular local rings are UFDs.

 A  resolution as in (2)
is said to be \it minimal \rm if $\alpha_i(F_i)\subseteq nF_{i-1}$. In this
case it is unique up to an isomorphism. The length of a minimal resolution
is the projective dimension of $M$, denoted pd$_S(M)$.

We say a Noetherian ring $R$ is regular if the localizations $R_P$ are
regular for evry prime ideal $P$ of $R$. 
Every polynomial ring or power series ring over a field is regular.

\medskip
\remark{Remark 2.2} Although regular rings were the building blocks
of commutative algebra, other types of rings were becoming
increasinly interesting and important. To explain these,
further definitions are helpful:

$\bullet$ An $m$-primary ideal $I$ is any ideal such that one of the
following equivalent conditions hold:
\roster
\item $R/I$ is dimension $0$;
\item $R/I$ is Artinian;
\item The nilradical of $I$ is $m$;
\item Supp$(R/I) = \{m\}$.
\endroster
\smallskip
$\bullet$ A regular sequence $x_1,...,x_t$ on an $R$-module $M$ is any set
of elements such that $(x_1,...,x_t)M\ne M$ and such that $x_1$ is not a
zero-divisor on $M$, $x_2$ is not a zero-divisor on the module $M/x_1M$, and
in general $x_i$ is not a zero-divisor on $M/(x_1,...,x_{i-1})M$.
\smallskip

$\bullet$ The depth of a finitely generated module over a local ring
$R$ is the maximal length of a regular sequence on the module.

$\bullet$ $M$ is said to be \it Cohen-Macaulay \rm if depth$(M) = \text{dim}(M)$.
(This is the largest possible value for the depth.)
\smallskip
$\bullet$ A s.o.p. (system of parameters) for a local ring $R$ of dimension $d$ is any $d
$
elements $x_1,...,x_d$ which generate an $m$-primary ideal.
\smallskip
$\bullet$ A local ring is Cohen-Macaulay iff every (equivalently one)
system of parameters forms a regular sequence.
\smallskip
$\bullet$ A local ring $R$ is regular iff the maximal ideal is
generated by a s.o.p.
\smallskip
$\bullet$ The height of a prime
ideal $P$ in a commutative ring $R$ is the supremum of integers $n$
such that there is a chain of distinct prime ideals
$P_0\subseteq P_1\subseteq ... \subseteq P_n = P$. The height of
an arbitrary ideal is the smallest height of any prime containing
the ideal.
\endremark

If $I$ has height $k$, then $I$ needs at least $k$ generators,
otherwise the dimension could not drop by $k$. This is the statement of
what is known as Krull's Height Theorem. Explicitly Krull's theorem
states:

\proclaim{Theorem 2.3} Let $S$ be a Noetherian ring and let $I$ be an ideal
generated by $k$ elements. Then every
minimal prime containing $I$ has height at most $k$. In particular
height$(I)\leq k$.
\endproclaim

The property of a ring being Cohen-Macaulay had been growing
in importance in the 1950s and early 1960s. Different groups
studied it under different names: semi-regular, Macaulay, even
Macaulay-Cohen. In particular, Northcott and Rees had been
studying the irreducibility of systems of parameters and relating
it to the Cohen-Macaulay property.

\smallskip
$\bullet$ An ideal $I$ is irreducible if $I\ne J\cap K$ for $I\subsetneq J,K$.
\smallskip
$\bullet$ Let $M$ be a module over a local ring $(R,m)$. The socle of a module
$M$, denoted by soc$(M)$, is the largest
subspace of $M$ whose $R$-module structure is that of a vector space, i.e., soc$(M) =$
 Ann$_Mm = \{y\in M|\, my = 0\,\}$.
\smallskip
\proclaim{Proposition 2.4} If $R$ is Artinian local, then $(0)$ is irreducible iff
soc$(R)$ is a 1-dimensional vector space.
\endproclaim

\demo{Proof} Let $V$ be the socle.
If dim$(V)\geq 2$ then simply pick two one dimensional
subspaces of $V$ which intersect in $0$ to show $0$ is not irreducible.

For the converse, the key point is that $R$ is Artinian implies
$V\subseteq R$ is essential. Any two ideals must intersect the
socle, but if it is one-dimensional they then contain the socle. \qed
\enddemo

For example, the ring $R = k[X,Y]/(X^2,Y^2)$ has a one-dimensional socle
generated by $XY$; every nonzero ideal in $R$ contains this element
and hence $0$ is irreducible. The intersection of any two nonzero ideals
must contain $XY$. On the other hand in $R = k[X,Y]/(X^2,XY,Y^2)$ the
socle is two dimensional generated by the images of $X$ and $Y$, and
$XR\cap YR = 0$. 

Irreducible ideals were already present in the famous proof of Emmy Noether
that ideals in Noetherian commutative rings have a primary decomposition \cite{No}. 
She first observed that the Noetherian property implied every ideal is a
finite intersection of irreducible ideals. (For example the ideal
$(X^2,XY,Y^2)$ above is the intersection of $(X,Y^2)$ and $(X^2,Y)$,
which are both irreducible.)  This `irreducible decomposition'
is closely related to the ideas in the theory of Gorenstein rings.

The results of Northcott and Rees are the following:

\proclaim{Theorem 2.5 \cite{NR} }
If $(R,m)$ is a local Noetherian ring such that the ideal
generated by an arbitrary system of parameters is irreducible,
then $R$ is Cohen-Macaulay.
\endproclaim

\proclaim{Theorem 2.6 \cite{NR}}
A regular local ring has the property that every ideal generated
by a system of parameters is irreducible.
\endproclaim

Combining these two theorems recovers a result first due to Cohen:
regular local rings are Cohen-Macaulay. In fact, as we shall see,
this irreducibility property characterizes Gorenstein rings.
>From this perspective, the theorems of Northcott and Rees say (among
other things) that regular rings are Gorenstein, and Gorenstein
rings are Cohen-Macaulay.

\bigskip
\head 3. Gorenstein Rings\endhead
\bigskip
Grothendieck (1957) \cite{Gr} developed duality theory for singular varieties.
He constructed a rank one module which on the non-singular
locus of an algebraic variety agreed with the differential forms of
degree $d$ (the dimension) of the variety.

Let $S$ be a polynomial ring  and $R$ be
a homomorphic image of $S$ of dimension $d$.
One can define a dualizing module (or \it canonical module\rm) by setting
$$
\omega_R = \text{Ext}_S^{n-d}(R,S)
$$
where the dimension of $S$ is $n$. When $R$ is non-singular, this
module is isomorphic to the module of differential forms of degree
$d$ on $R$.  The best duality is not obtained unless the
ring $R$ is homologically trivial in the sense that
$\text{Ext}_S^{j}(R,S) = 0$ for $j\ne n-d$. Such rings are exactly the
Cohen-Macaulay rings.

The duality is `perfect' on $R$-modules whose projective dimension
over the regular ring $S$ is smallest possible. (These are the
Cohen-Macaulay modules.) Let $M$ be such a module of dimension
$q$. Then $$\text{Ext}_R^i(M,\omega) = 0$$ for $i\ne d-q$,
$\text{Ext}_R^{d-q}(M,\omega)$ is again  Cohen-Macaulay and
$$ \text{Ext}_R^{d-q}(\text{Ext}_R^{d-q}(M,\omega),\omega))\cong M.$$

According to Bass in the ubiquity paper, 
Grothendieck introduced the following
definition:

\definition{Definition 3.1} Let $S$ be a regular local ring.
A local ring $R$ which is a homomorphic image of $S$ is \it Gorenstein \rm
if it is Cohen-Macaulay and its dualizing module (or \it canonical module\rm)
$\text{Ext}_S^{n-d}(R,S)$ is free (of rank 1), where $n = \text{dim}(S)$
and $d = \text{dim}(R)$. We say a  possibly non-local 
Noetherian ring $R$ (which we assume is the homomorphic image of a regular ring) is Gorenstein if $R_P$ is Gorenstein for all
prime ideals $P$.
\enddefinition
\smallskip
All plane curves are Gorenstein, but even more we shall see that 
all complete
intersections are Gorenstein. Definition 3.1 will be temporary in the
sense that the work of Bass made it clear how to extend the definition
to all local rings, and thereby all rings, without necessarily having
a canonical module.

The earliest printed reference I could find giving a definition of 
Gorenstein rings is in \cite{Se1, p. 2-11}. The date given for Serre's
seminar is November 21, 1960.  He defines an algebraic variety $W$ to
be Gorenstein if it is Cohen-Macaulay and the dualizing sheaf $\Omega_W$ 
is locally free of rank $1$. He goes on to say that for plane curves $R$
the equality dim$(T_Q/R_Q) = \text{dim}(R_Q/\frak C_Q)$ is equivalent
to being Gorenstein and cites \cite{Se3}. He further mentions that
complete intersections are Gorenstein in this context.  In \cite{Se3} (in 1959)
Serre had proved these claims and mentioned the papers of Rosenlicht, Samuel
and Gorenstein, as well as the work of Grothendieck on the duality approach.

With the definition of Grothendieck/Serre, it is not difficult to
see that complete intersections are Gorenstein.
We only need to use the free resolutions of such rings.
These resolutions give a powerful tool not only to prove results about
Gorenstein rings but to actually compute with them.

\proclaim{Proposition 3.2} Let $(S,m)$ be a regular local ring of dimension
$n$ and let $R = S/I$ be a quotient of $S$ of dimension $d$. Let
$$0\ra F_k\overset {\alpha_k}\to\ra F_{k-1}\overset {\alpha_{k-1}}\to\ra ...\ra F_0\ra R\ra 0$$
be a minimal free $S$-resolution of $R$. Then $R$ is Gorenstein iff
$k = n - d$, and $F_k\cong S$.
\endproclaim

\demo{Proof} We use a famous formula of Auslander and Buchsbaum to
begin the proof. If $S$ is a regular local ring and $M$ a finitely
generated module, then $$ \text{depth} (M) + \text{pd}_SM = \text{dim}(S)$$
Applying this formula with $R = M$ we obtain that $R$ is Cohen-Macaulay
iff depth$(R) = \text{dim}(R)$ (by definition) iff $k = n - d$. 

The module $\text{Ext}_S^{n-d}(R,S)$ can be computed from the
free $S$-resolution of $R$. We simply apply Hom$_S(\quad, S)$ to the
resolution and take homology. The $n-d$ homology is the cokernel of
the transpose of $\alpha_k$ from $F_{k-1}^*\ra F_k^*$, where
$(\quad)^* = \text{Hom}_S(\quad, S)$. If $R$ is Gorenstein, this is
free of rank $1$ by definition. The minimality of the resolution together
with Nakayama's lemma shows that the minimal number of generators of
$\text{Ext}_S^{n-d}(R,S)$ is precisely the rank of $F_k$. Hence it must be
rank $1$.

Conversely, suppose that the rank of $F_k$ is $1$. Since $I$ kills $R$,
it also kills $\text{Ext}_S^{n-d}(R,S)$. It follows that $\text{Ext}_S^{n-d}(R,S)$
is isomorphic to $S/J$ for some ideal $J\supseteq I$. However, the
vanishing of the other Ext groups gives that the transposed complex
$$0\ra F_0^*\ra ... \ra F_k^*$$
is actually acyclic, and hence a free $S$-resolution of $\text{Ext}_S^{n-d}(R,S)$.
Therefore $$\text{Ext}_S^{n-d}(\text{Ext}_S^{n-d}(R,S), S)\cong R$$  by
dualizing the complex back, and now the same reasoning shows that $J$ kills
$R$, i.e. $J\subseteq I$. Thus $\text{Ext}_S^{n-d}(R,S)\cong R$ and $R$ is
Gorenstein.\qed
\enddemo

Observe that the end of the proof proves that the resolution of a
Gorenstein quotient $R$ of $S$ is essentially self-dual; if you flip the
free resolution you again get a free resolution of $R$. In particular
the ranks of the free modules are symmetric around the `middle' of the
resolution.

\remark{Remark 3.3} Proposition 3.2 is especially easy to apply when the
dimension of $R$ is
$0$. In that case the depth of $R$ is $0$ which forces the
projective dimension to be $n$, so that the first condition of (3.2) is automatically
satisfied; one only needs to check that the rank of the last free
module in the minimal free resolution of $R$ is exactly $1$ to conclude
that $R$ is Gorenstein.
\endremark

\remark{Remark 3.4: The Graded Case}
If $R = \bigoplus_nR_n$ is graded with $R_0 = k$ a field,  and is
Noetherian, we
can write $R = S/I$ with $S = k[X_1,...,X_n]$ and $I$ is a homogeneous
ideal of $S$, where the $X_i$ are given weights. In this case the free resolution of $S/I$ over $S$
can be taken to be graded. If we twist\footnote{If $M$ is a graded
$S$-module, then the twist $M(t)$ is the same module but with the
different grading given by $M(t)_n = M_{t+n}$.} the free modules so that
the maps have degree $0$, then the resolution has the form:
$$
0\ra \sum_i S(-i)^{b_{ni}}\ra ... \ra \sum_i S(-i)^{b_{1i}}\ra S\ra S/I\ra 0.
$$

For example, $b_{1i}$ is the number of minimal generators
of $I$ in degree $i$. 

Let $M$ be the unique homogeneous maximal ideal
of $S$ generated by the $X_i$. Then
$R$ is Gorenstein iff $R_M$ is Gorenstein.  This follows from the
characterization above. One uses that the free resolution of
$R$ over $S$ can be taken to be graded, and hence the locus of primes
$P$ where
$R_P$ is Cohen-Macaulay is defined by the a homogeneous ideal
and the locus where the canonical module is free is also homogeneous
as the canonical module is homogeneous. This means that if there
exists a prime $P$ in $R$ such that $R_P$ is not Gorenstein, then
$P$ must contain either the homogeneous ideal defining the non Cohen-Macaulay
locus or the homogeneous ideal defining the locus where the canonical
module is not free. In either case, these loci are non-empty and will
also contain $M$, forcing $R_M$ not to be Gorenstein.
\endremark

To prove that complete intersections are Gorenstein, one only needs to
know about the Koszul complex: Let $x_1,...,x_n$ be a sequence of
elements in a ring $R$. The tensor product of the complexes
$$0\ra R\overset {x_i}\to\longrightarrow R\ra 0$$
is called the \it Koszul complex \rm \footnote{There is
another perhaps more fundamental way to think of the Koszul complex as
an exterior algebra. See \cite{BH, 1.6}.} of $x_1,...,x_n$.  This complex is a complex of
finitely generated free $R$-modules of length $n$. If the $x_i$ form a
regular sequence, then this complex provides a free resolution of
$R/(x_1,...,x_n)$.

\proclaim{Corollary 3.5} Let $S$ be a regular local ring, and let $I$ be
an ideal of height $k$ generated by $k$ elements. Then $R = S/I$ is
Gorenstein.
\endproclaim

\demo{Proof} Since $S$ is Cohen-Macaulay, $I$ is generated by a regular
sequence and the Koszul complex provides a free $S$-resolution of $R$.
The length of the resolution is $k$, and the last module in this
resolution is $S$. Applying Proposition 3.2  gives that $R$ is Gorenstein.\qed
\enddemo

The following sequence of implications is learned by all commutative
algebraists:
$$ \text{regular}\Rightarrow \text{complete intersection}\Rightarrow
\text{Gorenstein} \Rightarrow \text{Cohen-Macaulay}.$$ In general none
of these implications is reversible.  The last implication is
reversible if the ring is a UFD and a homomorphic image of a Gorenstein
ring, as we shall discuss later.

The language and results of free resolutions give a useful tool to
analyze Gorenstein rings, especially in low codimension.

\head 4. Examples and Low Codimension\endhead 
\medskip
Let $S$ be a regular local ring. If an ideal $I$ has height
one and $S/I$ is Gorenstein, then in particular $I$ is unmixed.
Since $S$ is a UFD, it follows that $I = (f)$ is principal, and
$S/I$ is a complete intersection.

 Let $I$ be a height two ideal
defining a Gorenstein quotient, i.e. $R = S/I$ is Gorenstein of
dimension $n-2$. The resolution for $R$ must look like:

$$0\ra S\ra S^l\ra S\ra S/I = R\ra 0.$$

But counting ranks shows that $l = 2$, which means that $I$ is
generated by $2$ elements and is a complete intersection. This
recovers a result of Serre \cite{Se1}.

A natural question is whether all Gorenstein rings are complete
intersections. The answer is a resounding NO, but it remains
an extremely important question how to tell if a given Gorenstein ring
is a complete intersection; this question plays a small but
important role in the work of Wiles, for example.

In codimension three, the resolution of a Gorenstein $R$ looks like:
$$0\ra S\ra S^l\ra S^l\ra S\ra S/I = R\ra 0.$$
The question is whether $l = 3$ is forced.
The answer is no. But surprisingly, the minimal number of generators
must be odd.

\example{Example 4.1} Let $S = k[X,Y,Z]$, and let $$I = 
(XY, XZ, YZ, X^2-Y^2, X^2-Z^2).$$ The nilradical of
$I$ is the maximal ideal $m = (X,Y,Z)$ which has height three.
But $I$ requires five generators and so is not a complete
intersection. On the other hand
$I$ \it is \rm Gorenstein. One can compute the resolution of $S/I$
over $S$ and see that the last Betti number is $1$. The resolution
looks like:
$$
0\rightarrow S\rightarrow S^5\rightarrow S^5\rightarrow S\rightarrow S/I\rightarrow 0.
$$

The graded resolution is:

$$
0\rightarrow S(-5)\rightarrow S^5(-3)\rightarrow S^5(-2)\rightarrow S\rightarrow S/I\rightarrow 0.
$$

\endexample

This example is a special case of a famous theorem of Buchsbaum and Eisenbud
\cite{BE}. To explain their statement we first recall what Pfaffians are.
Let $A$ be a skew-symmetric matrix of size $2n$ by $2n$. Then the determinant of
$A$ is the square of an element called the Pfaffian of $A$ (the sign is
determined by convention). If $A$ is skew symmetric and of size
$2n+1$, the ideal of Pfaffians of order $2i$ ($i\leq n$) of $A$ is the
ideal generated by the Pfaffians of the submatrices of $A$ obtained by
choosing $2i$ rows and the same $2i$ columns. The theorem of
Buchsbaum and Eisenbud states:

\proclaim{Theorem 4.2} Let $S$ be a regular local ring, and let $I$ be an ideal
in $S$ of height $3$. Set $R = S/I$. Then $R$ is Gorenstein iff $I$ is
generated by the 2n-order Pfaffians of a skew-symmetric $2n+1$ by $2n+1$
alternating matrix $A$.  In this case a minimal free resolution of $R$
over $S$ has the form,
$$0\ra S\ra S^{2n+1}\overset{A}\to\longrightarrow S^{2n+1}\ra S\ra R\ra 0.$$
\endproclaim

Unfortunately, there is no structure theorem for height four ideals
defining Gorenstein rings, although there has been a great deal of
work on this topic. See the references in \cite{VaVi} and \cite{KM}.

Thinking back to the first example of
monomial plane curves, it is natural to ask when the ring
$k[t^{n_1},...,t^{n_k}]$ is Gorenstein. 

\definition{Definition 4.3} A semigroup $T = <0,n_1,...,n_k>$ is said to
be \it symmetric \rm  if
there is a value $c\notin T$ such that $m\in T$ iff $c-m\notin T$.
\enddefinition

\proclaim{Theorem 4.4 \cite{HeK1}} The monomial curve $R = k[t^{n_1},...,t^{n_k}]$ is Gorenstein iff
$T = <0,n_1,...,n_k>$ is symmetric.
\endproclaim

It follows that every semigroup $T$ generated by two relatively prime
integers $a$ and $b$ is symmetric, since this is a plane curve,
hence Gorenstein.  For such examples, the value of $c$ as in (4.3) is $c = (a-1)(b-1) - 1$. In general
$c$ must be chosen one less than the least power of $t$ in the conductor.
This example can be used to give many examples of Gorenstein curves
which are not complete intersections or Pfaffians. The easiest such example is
$T = <0, 6, 7, 8, 9, 10>$. The ring $R = k[\![t^6,t^7,t^8,t^9,t^{10}]\!]$
is Gorenstein since $T$ is symmetric. The conductor is everything from $t^{12}$
and up, and one can take $c = 11$. The defining ideal when we represent $R$ as
a quotient of a polynomial ring in $5$ variables has height $4$ and
requires $9$ generators. The minimal resolution over the polynomial
ring $S$ has the form
$$
0\ra S\ra S^9\ra S^{16}\ra S^9\ra S\ra R\ra 0.
$$

See \cite{Br} for a discussion of symmetric semigroups.
For a non-Gorenstein example, we can take $R = k[t^3,t^4,t^5]$. The
semi-group is not symmetric. Writing $R$ as a quotient of a polynomial
ring in three variables, the defining ideal is a prime ideal of height
$2$ requiring $3$ generators. From the theorem of Serre  explained earlier,
one knows that this cannot be Gorenstein since it would have to be 
a complete intersection being height two. Alternatively, one could use
an amazing result due to Kunz \cite{Ku}, that almost complete intersections
are never Gorenstein. Almost complete intersection means minimally
generated by one more element than a complete intersection. There are
non-trivial examples of ideals generated by two more elements than
their height which define Gorenstein rings; the Pfaffian ideal of
$5$ generators gives such an examples, while whole classes were
constructed in \cite{HU2} (see also \cite{HM}).

Where does ubiquity fit into this? We need another historical thread
before the tapestry is complete. This last thread concerns the
theory of injective modules.

\bigskip
\head 5. Injective Modules and Matlis Duality\endhead
\bigskip

There is a `smallest' injective module containing an $R$-module
$M$, denoted $E_R(M)$, or just $E(M)$. It is called the \it injective
hull \rm of $M$. (It is the largest essential extension of $M$.)
Any module both esesential over $M$ and injective must be isomorphic
to the injective hull of $M$.
\medskip
\proclaim{Theorem 5.1 (Matlis, \cite{Ma})} Let $R$ be a Noetherian commutative
ring. Every injective module is a direct sum of indecomposable
injective modules, and the nonisomorphic indecomposable injective modules
are exactly (up to $\cong$) $E(R/P)$ for $P$ prime.
\endproclaim

Let $(R,m)$ be a $0$-dimensional Noetherian local ring. There is only
one indecomposable injective module, $E$, the injective hull of the
residue field of $R$. Matlis proved that the length of $E$, that is the
number of copies of the residue field $k = R/m$ in a filtration whose
quotients are $k$, is the same as that of $R$.
Matlis extended these ideas and came up with
what is now called Matlis duality.

\proclaim{Theorem 5.2} Let $(R,m)$ be a Noetherian local ring with
residue field $k$, $E = E_R(k)$, and completion $\widehat R$. Then there
is a $1-1$ arrow reversing correspondence from finitely generated $\widehat R$-modules
$M$ to Artinian $R$-modules $N$ given by $M\ra M^{\vee}$ and $N\ra N^{\vee}$
where $(\quad)^{\vee}:= \text{Hom}_R(\quad,E)$. Furthermore $M^{\vee\vee}\cong M$
and $N^{\vee\vee}\cong N$.  Restricting this correspondence to
the intersection of Artinian modules and finitely generated modules,
i.e. to modules of finite length, preserves length.
\endproclaim 

The work of Matlis allows one to understand the $0$-dimensional commutative
Noetherian rings which are self-injective. In fact, as the referee
pointed out to me, the equivalences (2) and (3) in Theorem 5.3 as well
as the basic duality that $0:0:I = I$ for ideals $I$ in a $0$-dimensional Gorenstein
ring, and further results on the lengths of ideals in such rings
were proved in a 1934 paper of W. Gr\" obner \cite{Gro}. He in turn
refers back to a 1913 paper of Macaulay \cite{Mac} giving similar results,
but says they are difficult to understand. Macaulay refers back to
the work of Lasker on primary decomposition.

\proclaim{Theorem 5.3} Let $(R,m)$ be a $0$-dimensional local 
ring with residue field $k$ and $E = E_R(k)$. The following are equivalent.
{\roster
\item $R$ is injective as an $R$-module.
\item $R\cong E$.
\item $(0)$ is an irreducible ideal.
\item The socle of $R$ is $1$-dimensional.
\item $R$ is Gorenstein in the sense of Grothendieck/Serre.
\endroster}
\endproclaim

\demo{Proof}
Clearly (2) implies (1). On the other hand, if $R$ is injective it
must be a sum of copies of $E$ by the theorem of Matlis. Since $R$
is indecomposable, it is then isomorphic to $E$.

Every $0$-dimensional
local ring is an essential extension of its socle $V$. For 
given any submodule of $R$, that is an ideal $I$ of $R$, there is a
least power $n$ of $m$ such that $m^nI = 0$. Then $m^{n-1}I\inc V\cap I$ and
is nonzero. Since the length of $R$ and $E$ agree, it then is clear
that $R$ is self-injective iff it's socle is $1$-dimensional, i.e.,
iff $0$ is irreducible. This shows the equivalence of (2), (3) and (4).

To see the equivalence of (5), first observe that $R$ is the homomorphic image
of a regular local ring $S$, using the Cohen structure theorem.
$R$ is Cohen-Macaulay because it is $0$ dimensional and of course has
depth $0$.  So it just needs to be shown that Ext$^n_S(R,S)$ is isomorphic
to $R$, where dim$(S) = n$. Equivalently, as shown in Proposition 3.2
we need to prove that the rank of the last free module in a minimal
free $S$ resolution of $R$ is exactly one. This rank is the dimension
of Tor$_S^n(k,R)$ by computing this Tor using the minimal free resolution of
$R$. However, $k$ has a free $S$-resolution by a Koszul complex since
the maximal ideal of $S$ is generated by a regular sequence as $S$ is
regular. Using the Koszul complex to compute this Tor yields the
isomorphism $$\text{Tor}_S^n(k,R)\cong \text{soc}(R).$$
Thus $R$ is Gorenstein iff the socle is $1$-dimensional. \qed
\enddemo
\smallskip
A natural question from Theorem 5.3 and from the second
theorem of Northcott and Rees is how to calculate the socle of a
system of parameters in a regular local ring $(R,m)$. 
For if one takes a regular local
ring and quotients by a system of parameters, then $0$ is irreducible
by Northcott and Rees, and the quotient ring is therefore a self-injective ring.
 It turns out there
is an excellent answer to how to calculate socles.

\proclaim{Theorem 5.4} Let $(R,m)$ be a regular local ring of dimension $d$,
and write $m = (x_1,...,x_d)$. Let $f_1,...,f_d$ be a 
system of parameters. Write $(x_1,...,x_d)A = (f_1,...,f_d)$ for some
$d$ by $d$ matrix $A$, and put $\Delta = \text{det}(A)$. Then the
image of $\Delta$ in $R/(f_1,...,f_d)$ generates the socle of this
Artinian algebra.
\endproclaim

Another way to state this theorem is $(f_1,...,f_d):_Rm = (f_1,...,f_d,\Delta)$.
Here $I:_RJ = \{r\in R|\, rJ\inc I\,\}$. Theorem 5.4 gives a very useful
computational device for computing socles. One easy example is to
compute a generator for the socle of $R/(x_1^{n_1},...,x_d^{n_d})$. In this
case the matrix $A$ can be taken to be a diagonal matrix whose $i$th entry
along the diagonal is $x_i^{n_i-1}$, and $\Delta = \prod_ix_i^{n_i-1}$.
Of course, when we lift the socle generator back to $R$ it is unique
only up to a unit multiple plus an element of the ideal we mod out.

In the case when $R$ is complete, regular, and contains a field, the Cohen
Structure Theorem gives that $R$ is isomorphic with a formal power series
over a field. In this case the theory of residues can be used to give
another description of the socle:

\proclaim{Theorem 5.5} Let $R= k[\![x_1,...,x_d]\!]$ be a formal power series
ring over a field $k$ of characteristic $0$, and let $f_1,...,f_d$
be a system of parameters
in $R$. Let $J$ be the Jacobian matrix whose $(i,j)$ entry is
$\frac {\partial f_i} {\partial x_j}$.
Set $\Delta = \text{det}(J)$. Then $\Delta$ generates the socle of
$R/(f_1,...,f_d)$.
\endproclaim

If the $f_i$ happen to be homogeneous polynomials, then this Theorem
follows directly from the previous one as Euler's formula can be
used to express the $f_i$ in terms of the $x_j$ using only units
times the rows of the Jacobian matrix.  However, in the nonhomogenous
case, this second theorem is extremely useful.

\remark{Remark 5.6} In the graded case in which $R = k[R_1]$ is 
$0$-dimensional and a homomorphic image of a polynomial ring $S =
k[X_1,...,X_n]$, one
can calculate the degree of the socle generator for $R$ by using
the free resolution of $R$ as an $S$-module. The resolution can be taken to
be of the form

$$
0\ra S(-m)\ra ... \ra \sum_i S(-i)^{b_{1i}}\ra S\ra R\ra 0.
$$

The last module will be a copy of $S$ twisted by an integer $m$ since
$R$ is Gorenstein.  The degree of the socle of $R$ is then 
exactly $m - n$. For example, in (4.1) we saw that the last twist
is $S(-5)$ in the resolution of the ring $R = k[X,Y,Z]/I$
where $I = (X^2-Y^2, X^2-Z^2, XY, XZ, YZ)$. Hence the socle of
$R$ sits in degree $ 2 = 5-3$. In fact the socle can be taken be
any quadric not in $I$.
\endremark
 
\bigskip
\head 6. Ubiquity \endhead
\bigskip
In the early 1960s, Bass had been studying properties of rings with
finite injective dimension. In particular, he related irreducibility
to finite injective dimension.

The following is from Bass \cite{Ba1} (1962): 
\medskip
\noindent\proclaim{Theorem 6.1} Let $R$ be a commutative Noetherian
local ring. The following conditions
are equivalent:
\roster
\item $R$ is Cohen-Macaulay and every system
of parameters generates an irreducible ideal.
\item \text{id}$_R(R) < \infty$.
\item id$_R(R) = \text{dim}(R)$.
\endroster
\endproclaim

One can see in retrospect how the themes of irreducibility, injective dimension
and complete intersections were coming together.  Most of the work
discussed above was available by the late 1950s. The Auslander-Buchsbaum
work was done in 1957-58, Grothendieck's work on duality in 1957, Serre's
work in 1957-1960, Matlis's
work on injective modules in 1957, and the work of Northcott and Rees on
irreducible systems of parameters was done in 1957. The work of
Rosenlicht and Gorenstein was done about 1952. According to Bass, 
Serre pointed out to him that the rings of finite injective dimension
were, at least in the geometric context, simply the Gorenstein rings of
Grothendieck.  To quote Bass's paper he remarks, `Gorenstein rings, it is now
clear, have enjoyed such a variety of manifestations as to justify, perhaps,
a survey of their relevance to various situations and problems.' The rest,
as one says, is history.

Bass put all this together in the ubiquity paper in 1963 \cite{Ba2}. It remains
one of the most read papers in commutative algebra.

\proclaim{Theorem 6.2 Ubiquity} Let $(R,m)$ be a Noetherian local ring. The following
are equivalent.
{\roster
\item If $R$ is the homomorphic image of a regular local ring, then
$R$ is Gorenstein in the sense of Grothendieck/Serre.
\item id$_R(R) < \infty$.
\item id$_R(R) = \text{dim}(R)$.
\item $R$ is Cohen-Macaulay and some system of parameters generates an
irreducible ideal.
\item $R$ is Cohen-Macaulay and every system of parameters generates an
irreducible ideal.
\item If $0\ra R\ra E^0\ra ...\ra E^h ...\ra...$
is a minimal injective resolution of $R$, then for each $h\geq 0$,
$E^h\cong \sum_{\text{height}(p) = h} E_R(R/p).$
\endroster}
\endproclaim

By Theorem 2.5 of Northcott and Rees one can remove the assumption that
$R$ be Cohen-Macaulay in (5). 
We can now take any of these equivalent properties to be the definition
of Gorenstein in the local case. The usual one chosen is the second,
that the ring have finite injective dimension over itself. In the
non-local case we say $R$ is Gorenstein if $R_P$ is Gorenstein for
every prime $P$ in $R$. If $R$ has finite Krull dimension this
is equivalent to saying $R$ has finite injective dimension over itself.

\example{Example 6.3} Let $(R,m)$ be a one dimensional Gorenstein local ring.
The equivalent conditions above say that the injective resolution of
$R$ looks like $$0\ra R\ra K\ra E\ra 0$$
where $K$ is the fraction field of $R$ (this is the injective hull of $R$)
and $E$ is the injective hull of the residue field $k = R/m$. This gives
a nice description of $E$ as $K/R$. If we further assume that $R = k[t]_{(t)}$,
a regular local ring, then one can identify $K/R$ with the inverse powers
of $t$. This module is an essential extension of the residue field $k$
which sits in $K/R$ as the $R$-span of $t^{-1}$. This example generalizes
to higher dimensions in the form of the so-called `inverse systems' of
Macaulay. See section $8$.
\endexample

Moreover in dimension one for Noetherian domains 
with finite integral closure,
Bass was able to generalize the results of Gorenstein and Rosenlicht to
arbitrary Gorenstein rings:

\proclaim{Theorem 6.4} Let $S$ be a regular local ring and let $R = S/I$
be a $1$-dimensional Gorenstein domain. Assume that the integral closure of
$R$ is a finite $R$-module $($this is automatic for rings essentially of
finite type over a field, or for complete local rings$)$.
 Let $T$ be the integral
closure of $R$ and let $\frak C$ be the conductor. Then
$$ \lambda(R/\frak C) = \lambda(T/R).$$
\endproclaim

\demo{Proof} Let $\W$ be the canonical module of $R$, namely
Ext$_S^{n-1}(R,S)$ where $n = \text{dim}(S)$. Consider the exact
sequence, $0\ra R\ra T\ra D\ra 0$, where $D\cong T/R$.
Applying $(\quad)^* = \text{Hom}_R(\quad, R)$ we get the sequence,
$$
D^*\ra T^*\ra R^*\ra \text{Ext}^1_R(D,R)\ra \text{Ext}^1_R(T,R).
$$

Since $T$ is a Cohen-Macaulay of dimension $1$, the last Ext vanishes.
As $D$ is torsion $D^* = 0$. Moreover $R^*\cong R$ and $T^*\cong \frak C$,
the isomorphism given by sending a homomorphism to its evaluation at
a fixed nonzero element. It follows that $R/\frak C\cong \text{Ext}^1_R(D,R)$,
and it only remains to prove that the length of $\text{Ext}^1_R(D,R)$ is
the same as the length of $D$. This follows since $R$ is Gorenstein. We
can use the injective resolution of $R$ as in Example 6.3 to find that
$$\text{Ext}^1_R(D,R)\cong \text{Hom}_R(D,E).$$
By Matlis duality the length of this last module is the same as the
length of $D$. \qed 

\enddemo

\remark{Remark} In fact the converse to Theorem 6.4 is also true.
The converse appears in \cite{HeK1}, Corollary 3.7. The converse may have been known
to Serre, but this is as far as I know the first place the proof appears in print,
and is due to Herzog.
\endremark

\head 7. Homological Themes\endhead
\bigskip
One of the remarks Bass makes in the ubiquity paper is the following:
`It seems conceivable that, say for $A$ local, there exist finitely
generated $M\ne 0$ with finite injective dimension only if $A$ is
Cohen-Macaulay. The converse is true for if ($B = A$ modulo a system
of parameters) and $M$ is a finitely generated non-zero $B$-injective
module then inj dim$_A(M) < \infty$.'

This became known as Bass's conjecture, one of a celebrated group of
`homological conjectures'. The search for solutions to these problems
was a driving force behind commutative algebra in the late 1960s and
1970s, and continues even today. However, Bass's conjecture 
has now been solved positively, by Peskine and
Szpiro \cite{PS1} in the geometric case, by Hochster \cite{Ho1} for all local
rings containing a field, and by Paul Roberts \cite{Rob1} in mixed
characteristic.

The following is due to Peskine and Szpiro \cite{PS1, (5.5) and (5.7)}
based on their work on the conjecture of Bass:

\proclaim{Theorem 7.1}
Let $(R,m)$ be a local ring. The following are equivalent:
{\roster
\item R is Gorenstein.
\item $\exists$ an ideal $I$ such that id$_R(R/I) < \infty$.
\item $\exists$ an $m$-primary ideal $I$ such that pd$_R(R/I) < \infty$ and
$I$ is irreducible.
\endroster}
\endproclaim

Other homological themes were introduced in the ubiquity paper which are
now part of the standard landscape. A module $M$ of a local ring $R$
 is said to be a $k$th syzygy
if there is an exact sequence,
$$0\ra M\ra F_{k-1}\ra ... \ra F_0$$
where the $F_i$ are finitely generated free $R$-modules.
It is important to know intrinsically when a given module is a $k$th
syzygy. For Gorenstein local rings there is a good answer 
due to Bass \cite{Ba2, Theorem 8.2}:

\proclaim{Theorem 7.2} Let $R$ be a Gorenstein local ring, and $k\geq 2$.
A finitely generated $R$-module $M$ is a $k$th syzygy iff $M$ is reflexive
and Ext$^i_R(M^*,R) = 0$ for all $i\leq k - 1$.
\endproclaim

The idea of the proof is to take a resolution of the dual module
$M^* : = \text{Hom}_R(M,R)$ and then dualize it, using that $M\cong M^{**}$
since $M$ is reflexive. A modern version of this result was given by
Auslander and Bridger \cite{AB}, see for example 
\cite{EG, Theorem 3.8}.

Other important homological objects introduced by Bass in ubiquity were
the Bass numbers, $\mu^i(p,M) = \text{dim}_{k(p)}\text{Ext}^i_{R_p}(k(p),M_p)$,
where $k(p) = R_p/pR_p$ are the residue fields of the localizations at
prime ideals of $R$. These numbers play a very important role in
understanding the injective resolution of $M$, and are the subject of
much work. For example, another equivalent condition for a local Noetherian
ring $(R,m)$ to be Gorenstein is that $R$ be Cohen-Macaulay and the $d$th Bass number,
$\mu^d(m,R) = 1$. This was proved by Bass. Vasconcelos conjectured that
one could delete the hypothesis that $R$ be Cohen-Macaulay. This was proved by
Paul Roberts in 1983 \cite{Rob2}. Another line of work inspired by
the work of Bass is that of  Avramov and Foxby, who have developed and studied
the concept of homomorphisms between rings being Gorenstein. 
For example see \cite{AvF}
and the references there. 

\bigskip
\centerline{\bf 8. Inverse Powers and $0$-dimensional Gorenstein Rings }
\bigskip

The study of Gorenstein rings can be approached by first trying to
understand the $0$-dimensional Gorenstein rings. These are all the
$0$-dimensional Artinian rings which are injective as modules over
themselves. Equivalently they are exactly the $0$-dimensional Artinian rings
with a $1$-dimensional socle.  Given any local Gorenstein ring of
arbitrary dimension, one can always mod out the ideal generated by
a system of parameters and obtain a $0$-dimensional Gorenstein ring.
In the $0$-dimensional case, one has the formidable results of Matlis
to help. 

\example{Example 8.1} Let $R = k[X_1,...,X_n]$, or a power series ring
$k[\![X_1,...,X_n]\!]$. An injective hull of the
residue field $k = R/(X_1,...,X_n)$ is given by the \it inverse
powers\rm, $$E = k[X_1^{-1},...,X_n^{-1}]$$ with the action
being determined by
$$
X_i(X_1^{-a_1}\cdots X_n^{-a_n}) = X_1^{-a_1}\cdots X_i^{-a_i+1}\cdots  X_n^{-
a_n}
$$
if $a_i\geq 1$. If $a_i = 0$, then the product is $0$. Notice the
copy of $k$ inside $E$ is exactly the span of $1$.

The history of inverse powers is a long one, but their inception
as far as this author knows is in work of Macaulay. He was indeed
far-sighted!
\endexample

\example{Example 8.2} To obtain the injective hull of the residue field of
 a graded quotient
$$k[X_1,...,X_n]/I$$ one simply takes Hom$_S(S/I, E)$. This follows easily
from the fact that this module is injective over $S/I$ (using Hom-tensor
adjointness) and is also essential over $k$, since it can be identified
with a submodule of $E$ which is already essential over $k$.
This is naturally identified with the elements of $E$ killed by $I$.
One sees that the graded structure of $E$ is simply that of $R$
upside down.
\endexample

$R$ is $0$-dimensional Gorenstein local iff $R\cong E(k)$.
In this case, $R\cong \text{Hom}_S(R,E)$ which means that the annihinlator of
$R$ in the inverse powers is cyclic. The converse also holds. The
$m$-primary ideals $I$ such that $R = S/I$ are Gorenstein are
EXACTLY the annihilators of single elements in the inverse powers.
(See Proposition 8.4 below.)

\example{Example 8.3}
Let $F = X^{-2}+Y^{-2}+Z^{-2}$. The set of
elements in $S$ killed by $F$ is the ideal
$$I = (X^2-Y^2, Y^2-Z^2, XY, XZ, YZ).$$
\endexample

This is the height $3$ ideal defining the Gorenstein ring that we considered in Example
4.1. 

It is not difficult to prove that if $F$ is a quadratic form in
$X^{-1},Y^{-1},Z^{-1}$, then the rank of $F$ as a quadratic form
is $3$ iff the corresponding ideal has $5$ generators, while
if the rank is at most $2$, the corresponding ideal is a complete
intersection generated by $3$ polynomials. See \cite{Ei, p. 551, exer. 21.6}.

In characteristic $0$, there is an important form of
inverse powers using the ring of differential operators. Let $k$ be
a field of characteristic $0$ and let $T = k[y_1,...,y_n]$ be a
polynomial ring over $k$. A polynomial differential operator with
constant coefficients is an operator on $T$ of the form
$$
D = \sum_{i_1,...,i_r} a_{i_1,...,i_r}(\frac {\partial} {\partial y_1})^{i_1}\cdots (\frac {\partial} {\partial y_n})^{i_n}$$
with $a_{i_1,...,i_r}\in k$. We think of $D$ as an element in a new
polynomial ring $S = k[x_1,...,x_n]$ acting on $T$ by letting $x_i$
act as $\frac {\partial} {\partial y_i}$. Since the partials commute
with each other this makes sense. Fix $f\in T$, and let $I = I_f$ be
the ideal of all elements $D\in S$ such that $Df = 0$. Then $I$ is an
ideal in $S$ primary to $(x_1,...,x_n)$ such that $S/I$ is $0$-dimensional
Gorenstein. Furthermore, all such $0$-dimensional Gorenstein quotients
of $S$ arise in this manner. The point is that we can identify $T$
with the injective hull of the residue field of $S$. See
\cite{Ei, exer. 21.7}. This point of view has been very important for
work on the Hilbert scheme of $0$-dimensional Gorenstein quotients of
fixed dimension. See for example the paper of Iarrobino \cite{Ia} and
its references.

There is a related way to construct all $0$-dimensional Gorenstein
local rings containing a field. Start with a power series ring
$S = k[\![X_1,...,X_n]\!]$. Choose an arbitrary system of parameters
$f_1,...,f_n$ in $S$ and an arbitrary element $g\notin (f_1,...,f_n)$.
Set $I = (f_1,...,f_n):_Sg$. Then $R = S/I$ is a $0$-dimensional
Gorenstein ring. Conversely, all $0$-dimensional Gorenstein rings
are of this form.
The proof follows at once from the following proposition:

\proclaim{Proposition 8.4} Let $R$ be a $0$-dimensional Gorenstein ring,
and let $J\inc R$. Then $R/J$ is Gorenstein iff $J = 0:Rg$ for some
nonzero element $g\in R$. Precisely, $Rg = 0:_RJ$.
\endproclaim

\demo{Proof} Since $R$ is Gorenstein, it is isomorphic to its own
injective hull. It follows that the injective hull of $R/J$ is 
$\text{Hom}_R(R/J,R)\cong 0:_RJ$, and $R/J$ is Gorenstein iff
$0:_RJ\cong R/J$ iff $0:_RJ = gR$ and $0:_Rg = J$. \qed
\enddemo

For example, to construct the Gorenstein quotient $k[X,Y,Z]/ I$ where
$I = (X^2-Y^2, Y^2-Z^2, XY, XZ, YZ)$, we can do the following.
First take the complete intersection
quotient $T = k[X,Y,Z]/(X^3,Y^3,Z^3)$. Set $g = X^2Y^2+X^2Z^2+Y^2Z^2$. Then
$I = (X^3,Y^3,Z^3):g$ and is therefore Gorenstein.  ALL Gorenstein
quotients $R/K$ such that $(X^3,Y^3,Z^3)\inc K$ arise in this fashion:
$K = (X^3,Y^3,Z^3):f$ for some $f\in k[X,Y,Z]$ and $f$ is unique up to
an element of $(X^3,Y^3,Z^3)$ and a unit multiple.

Understanding $0$-dimensional Gorenstein rings is extremely important
if one uses Gorenstein rings as a tool. For a great many problems one
can reduce to a $0$-dimensional Gorenstein ring. A natural question in
this regard is the following: let $(R,m)$ be a Noetherian local ring.
When does there exist a descending sequence of $m$-primary ideals $I_n$
such that $R/I_n$ are Gorenstein for all $n$ (equivalently $I_n$ are
irreducible) and such that $I_n\inc m^n$? Melvin Hochster completely
answered this question. He called such rings \it approximately Gorenstein\rm.

\proclaim{Theorem 8.5 \cite{Ho3, 1.2, 1.6}} A local 
Noetherian ring $(R,m)$ is approximately Gorenstein if
and only if
its $m$-adic completion is approximately Gorenstein. 
An excellent local ring $R$
$($so, for example, a complete local Noetherian ring$)$ with
dim$(R)\geq 1$ is  approximately Gorenstein if and only if the
following two conditions hold:
\roster
\item  depth($R)\geq 1$.
\item If $P\in \text{Ass}(R)$ and dim$(R/P) = 1$, then $R/P\oplus R/P$
is not embeddable in $R$.
\endroster
\endproclaim

\bigskip
\head 9. Hilbert Functions\endhead
\bigskip

The \it Hilbert function \rm  of a graded ring over a field is the function
$H(n) =  \text{dim}_kR_n$. The generating function for the Hilbert
function is called the Hilbert series, i.e.
the series
$$
F(R,t) = \sum_{n\geq 0}(\text{dim}_kR_n)t^n
$$
  If $R$ is generated by one-forms and is
Noetherian, we can always write
$$
F(R,t) = \frac {h_0 + h_1t + ... +h_lt^l} {(1-t)^d}
$$
where $d = \text{dim}R$, the Krull dimension of $R$.
The vector of integers $(h_0,h_1,...,h_l)$ is called the \it h-vector \rm
of $R$. For example, if $R$ is $0$-dimensional, then $h_i = H(i)$.

Let $R = k[R_1]$. We can write $R = S/I$ where $S = k[X_1,...,X_n]$ and $I$ is a homogeneous
ideal of $S$. In this case the free resolution of $S/I$ over $S$
can be taken to be graded, as in Section 3 above.
The resolution has the form
$$
0\ra \sum_i S(-i)^{b_{ni}}\ra ... \ra \sum_i S(-i)^{b_{1i}}\ra S\ra S/I\ra 0.
$$

Fixing a degree $n$ allows one to compute the
dimension of $(S/I)_n$ as the alternating sum $\sum_{j,i} (-1)^j(S(-i)^{b_{ji}})_n$. Each term in this sum is easily computatable as the dimension of
the
forms of a certain degree in the polynomial ring.

To do this calculation, it is more convenient to write
$$ F_i = S(-c_{1i})\oplus ... \oplus S(-c_{n_ii}).
$$
Then the Hilbert series for $S/I$ is exactly
$$
\frac {\sum_{i = 0}^t(-1)^i(t^{c_{1i}} + ... + t^{_{c_{n_ii}}})} {(1-t)^n}.
$$

The free resolution of a Gorenstein ring is essentially symmetric.
The dual of the resolution of a Cohen-Macaulay ring always gives a free
resolution of Ext$^{n-d}_S(R,S)$. 
When $R$ is further assumed to be Gorenstein the canonical module is
isomorphic to a twist $R(-p)$ of $R$, and the flipped resolution is isomorphic to the
orginial resolution, with appropriate shifts. This basic duality should make the
following theorem no surprise:

\proclaim{Theorem 9.1} A $0$-dimensional graded Gorenstein
ring $R = R_0\oplus R_1\oplus ... \oplus R_t$ $($with $k = R_0$ a field and
$R_t\ne 0)$
has symmetric Hilbert function, i.e. $H(i) = H(t-i)$ for all $i = 0,...,t$.
\endproclaim

One sees this by considering the pairing $R_i \times R_{t-i}\ra R_t = k$.
Since $R$ is essential over the socle, it is not difficult to show
that this pairing is perfect and identifies $R_i$ with the dual of
$R_{t-i}$.

\example{Example 9.2}
Let us compute the Hilbert function of the Example 4.1. The ideal
is defined by five quadrics in three variables, and 
every cubic is in the defining ideal. Hence the Hilbert
function is $(1,3,1)$.
\endexample

\example{Example 9.3}
Let $f_1,...,f_n$ be a homogeneous system of parameters
in $S = k[X_1,...,X_n]$ of degrees
$m_1,...,m_n$ respectively generating an ideal $I$. The Hilbert function can be computed
from the graded free resolution of $S/I$, and this is just
the Koszul complex. The various graded free modules in this
resolution depend only upon the degrees of the $f_i$, and so the Hilbert function is
the same as that of $X_1^{m_1},...,X_n^{m_n}$. This is easily computed as above.
The socle is generated in degree $\sum_i(m_i-1)$.
\endexample

This symmetry and the fact the socle is one-dimensional is behind
the classical applications of the Gorenstein property outlined in
\cite{EGH}. That article details nine versions of the so-called
Cayley-Bacharach Theorem, beginning with the famous theorem of
Pappus proved in the fourth century A.D., and ending with a common
generalization of all of them: that polynomial rings are Gorenstein!
It is worth repeating one of these avatars here and hopefully
hooking the reader to look at \cite{EGH}. The following theorem was
proved by Chasles in 1885:

\proclaim{Theorem 9.4} Let $X_1, X_2\inc \bold P^2$ be cubic plane curves
meeting in exactly nine points. If $X\inc \bold P^2$ is any cubic containing
eight of these points, then it contains the ninth as well.
\endproclaim

\demo{Proof}
First note that $X_1$ and $X_2$ meet in the maximum number of
points by B\'ezout's theorem. We interpret the theorem 
algebraically. Let $S = k[X,Y,Z]$ be the homogeneous coordinate
ring of the projective plane. $X_1$ and $X_2$ correspond to 
the vanishing loci of two homogeneous cubics $F_1, F_2$. The ideal
$I = (F_1, F_2) = Q_1\cap ... \cap Q_9$ where $Q_i$ is defined
by two linear forms. These correspond to the nine points in which
$X_1$ and $X_2$ meet.  Assume that another cubic, $G$, is contained
in $Q_1\cap ... \cap Q_8$. We want to prove that $G\in Q_9$, i.e.
that $G\in I$. 

We proceed to cut down to Artinian quotients by killing a general
linear form. Set $T = S/(\ell)$, where $\ell$ is a general linear
form. Of course $T$ is simply a polynomial ring in two variables.
We write $f_1,f_2,g$ for the images of $F_1, F_2, G$ in this
new ring. $R = T/(f_1, f_2)$ is a complete intersection whose Hilbert
function is $(1,2,3,2,1)$.
The socle of $T/(f_1, f_2)$ sits in degree $4$. If $G\notin I$,
then $g\notin (f_1, f_2)$, and $g$ sits in degree $3$. Since 
$R$ is Gorenstein, $Rg$ must contain the socle. As $g$ sits in degree
three, the Hilbert function of the ideal $Rg$ is $(0,0,0,1,1)$ and hence $R/Rg$ has
$k$-dimension  seven. The Hilbert function of
$R/Rg$ must be $(1,2,3,1)$. However,
by assumption $(F_1, F_2, G)$ is contained in the eight linear
ideals $Q_1,...,Q_8$. When we cut by a general linear form, it follows
that the dimension must be \it at least \rm $8$! This contradiction
proves that $G\in I$. \qed
\enddemo

A beautiful result of Richard Stanley shows that the symmetry of
the Hilbert function is not only necessary but even a sufficient condition 
for a Cohen-Macaulay graded domain to be Gorenstein.

\proclaim{Theorem 9.5 \cite{St3}} If $R = k[R_1]$ is a Cohen-Macaulay
domain of dimension $d$, then $R$ is Gorenstein iff
$$
F(R, t) = (-1)^dt^lF(R, 1/t)
$$
for some $l\in \bold Z$.
\endproclaim

The functional equation in Theorem 9.5 is equivalent to the symmetry of
the $h$-vector.
\smallskip
\example{Example 9.6} The converse is NOT true if the ring is not a domain.
For instance, $I = (X^3,XY,Y^2)$ has Hilbert function $(1,2,1)$ and satisfies
the equation $F(t) = t^2F(1/t)$, but is not Gorenstein.
\endexample

It is an open problem posed by Stanley to characterize which sequences 
of integers can be the h-vector of a Gorenstein graded algebra.
Symmetry is a necessary but not sufficient condition. Stanley gives
an example of a Gorenstein $0$-dimensional graded ring with
Hilbert function $(1,13,12,13,1)$.  The lack of unimodality is the
interesting point in this example.

Stanley's theorem was at least partially motivated, as far as I can tell,
by the uses he makes of it. One of them is to prove that certain invariants of
tori are Gorenstein. Invariants of group actions provide a rich source of
Gorenstein rings, as we will discuss in the next section.
\bigskip

\head 10. Invariants and Gorenstein rings \endhead
\bigskip

There are a great many sources of Gorenstein rings. Complete interesections
are the most relevant. After all every finitely generated $k$-algebra
which is a domain is birationally a complete intersection, hence birationally
Gorenstein. The Gorenstein property behaves well under flat maps, and
often fibers of such maps are Gorenstein. For example, if $K$ is a finitely
generated field extension of a field $k$ and $L$ is an arbitrary field
extension of $k$, then $K\otimes_kL$ is Gorenstein. 

A huge  source of Gorenstein rings comes from invariants of groups.  We fix notation and
definitions. Let $G$ be a closed algebraic subgroup of the general linear group
$GL(V)$, where $V$ is a finite dimensional vector space over a ground
field $k$. We obtain a natural action of $G$ on the polynomial ring $k[V] = S$,
and we denote the ring of invariant polynomials by $R = S^G$. Thus,
$$
S^G = \{f\in k[V]\,|\, g(f) = f \, \text{for all}\, g\in G\}.
$$ 

We will focus on groups $G$ which are linearly reductive. This means that every
$G$-module $V$ is a direct sum of simple modules. Equivalently, every $G$-submodule
of $V$ has a $G$-stable complement. Examples of linearly reductive groups
include finite groups whose order is invertible, tori, the classical groups
in characteristic $0$, and semisimple groups in characteristic $0$.  When
$G$ is linearly reductive, there is a retraction $S\rightarrow S^G$ called
the \it Reynolds operator \rm which is $S^G$-linear. It follows that the
ring of invariants $R$ must be Noetherian.  For a chain of ideals in $R$,
when extended to $S$ will stabilize, and then applying the retraction
stabilizes the chain in $R$. But even more is true: a famous theorem
of Hochster and Roberts \cite{HoR} says that the ring of invariants is
even Cohen-Macaulay.

\proclaim{Theorem 10.1} Let $S$ be a Noetherian $k$-algebra which is regular
$($e.g. a polynomial ring over $k)$ and let $G$ be a linearly reductive
linear algebraic group acting $k$-rationally on $S$. Then $S^G$ is
Cohen-Macaulay.
\endproclaim

A natural question is to ask when $S^G$ is Gorenstein. The following
question was the focus of several authors (see \cite{Ho2, St1, St2, Wa} and
their references).

\remark{Question 10.2} Let $k$ be an algebraically closed field and let
$G$ be a linearly reductive algebraic group over $k$ acting linearly
on a polynomial ring $S$ over $k$ such that the det is the trivial
character. Is $S^G$ Gorenstein?
\endremark

Thus, for example, when $G$ is not only in $GL(V)$, but in $SL(V)$, 
the ring of invariants should be Gorenstein. This question was motivated
because of the known examples: it was
proved for finite groups whose order is invertible by Watanabe \cite{Wa},
for an algebraic torus by Stanley \cite{St1}, or for connected semisimple
groups by a result of Murthy \cite{Mu} that Cohen-Macaulay rings which
are UFDs are Gorenstein, provided they have a canonical module.
However, this question has a negative solution.  
Friedrich Knop in 1989 \cite{Kn} gave a counterexample to Question 10.2
and has basically completed described when the ring of invariants is
Gorenstein.

Murthy's result quoted in the above paragraph is straightforward
based on the theory of the canonical module. The canonical module of an
integrally closed Noetherian ring $R$ is
a rank $1$  reflexive module which therefore lives in the class group.
If $R$ is a UFD, the class group is trivial and so is the canonical
module. But this means that it is a free $R$ module, which together
with $R$ being Cohen-Macaulay gives Gorenstein.\footnote{ A partial converse
to Murthy's result was given by Ulrich \cite{U}. His construction
uses the theory of liaison, or linkage, in which ideals defining Gorenstein rings
play an extremely important role. See for example, \cite{PS2} or
\cite{HU1}.}  Griffith has shown in \cite{Gr, Theorem 2.1}
that if $S$ is a local ring which is a UFD, which is finite over
a regular local ring $R$ in such a way that the extension of fraction
fields is Galois, then $S$ is necessarily a complete intersection.
See also \cite{AvB}.

\bigskip
\head 11. Ubiquity and Module Theory\endhead
\bigskip

In trying to recapture the material that led Bass to write his
ubiquity paper it seems that his choice of the word `ubiquity'
came from that fact that the Gorenstein property was
arising in totally different contexts from many authors; from
Northcott and Rees' work on irreducible systems of parameters and
Cohen-Macaulay rings, from the work of Gorenstein and Rosenlicht
on plane curves and complete intersections, from the work of
Grothendieck and Serre on duality, and from his own work on
rings of finite injective dimension. They were indeed ubiquitous.

But there are other ways in which they are ubiquitous. Every complete local
domain or finitely generated domain over a field is birationally a
complete intersection; hence up to integral closure every such ring is Gorenstein.

Every Cohen-Macaulay ring $R$ with a canonical module $\omega$ is
actually a Gorenstein ring up to nilpotents by using the idealization
idea of Nagata.  Specifically,
 we form a new ring consisting of 2 by 2 matrices whose
lower left hand corner is $0$, whose diagonal is a constant element
of $R$, and whose upper left component is an arbitrary element of
$\omega_R$. Pictorially, elements look like
\smallskip
$$
\left(\matrix
r & y\\
0 & r
\endmatrix\right)
$$
where $r\in R$ and $y\in \omega_R$.
It is easy to see this ring $S$ is commutative, $\omega$ is an ideal
of $S$, $S/\omega\cong R$, and $\omega^2 = 0$. What is amazing is that
$S$ is Gorenstein! Up to `radical' this transfers questions about
Cohen-Macaulay rings to questions about Gorenstein rings.

Even if a local ring $R$ is not Cohen-Macaulay, we can still closely
approximate it by a Gorenstein ring if it is the homomorphic image of
a regular local ring $S$ (e.g. in the geometric case or the complete
case).  Write $R = S/I$ and choose a maximal regular sequence
$f_1,...,f_g\in I$ where $g = \text{height}(I)$. Then the ring
$T = S/(f_1,...,f_g)$ is Gorenstein, being a complete intersection,
$T$ maps onto $R$, and dim$(T) = \text{dim}(R)$. 

Regular rings are the most basic rings in the study of commutative rings.
However, Gorenstein rings are the next most basic, and as the examples
above demonstrate, one can approximate arbitrary local commutative
rings quite closely by Gorenstein rings. Moreover, the module theory
over Gorenstein rings is as close to that of regular rings as one 
might hope.

Recall that finitely generated modules over a regular local ring are
of finite projective dimension. This characterizes regular local rings,
so we cannot hope to achieve this over non-regular rings. However there
are two `approximation' theorems for finitely generated modules over
Gorenstein local rings which show that up to Cohen-Macaulay modules,
we can still approximate such modules by modules of finite projective
dimension. The first theorem goes back to Auslander and Bridger \cite{AB}.

\proclaim{Theorem 11.1} Let $(R,m)$ be a Gorenstein local ring, and let $M$
be a finitely generated $R$-module of dimension equal to the
dimension of $R$. Then there exists an exact sequence,
$$
0\ra C\ra M\oplus F\ra Q\ra 0$$
where $F$ is finitely generated and free, $C$ is a Cohen-Macaulay
module of maximal dimension, and $Q$ is a module of finite projective
dimension.
\endproclaim

A second and similar theorem in spirit is due to Auslander and Buchweitz \cite{ABu}:

\proclaim{Theorem 11.2} Let $(R,m)$ be a Gorenstein local ring, and let $M$
be a finitely generated $R$-module. Then there is an exact sequence
$$
0\ra Q\ra C\ra M\ra 0$$
such that $C$ is a finitely generated Cohen-Macaulay module of maximal dimension and $Q$
has finite projective dimension.
\endproclaim

To put this in context, any Cohen-Macaulay module of maximal dimension over
a regular local ring is free, so both of these results recover that
over a regular local ring finitely generated modules have finite
projective dimension. The point is that over a Gorenstein ring, the
study of modules often reduces to the study of modules of finite
projective dimension and Cohen-Macaulay modules of maximal dimension.
This is certainly the best one can hope for.

The modules $C$ in the statement of (11.2) are called Cohen-Macaulay approximations of
$M$ and are a topic of much current interest. Theorem 11.2 can be thought
of as a generalization of an argument known as Serre's trick, which he
used in the study of projective modules. 
Let $M$ be a  finitely generated
module over a ring $R$ and choose generators for Ext$^1_R(M,R)$. If there
are $n$ such generators, then one can use the Yoneda definition of Ext to
create an exact sequence $0\ra R^n\ra N\ra M\ra 0$ where Ext$^1_R(N,R) = 0$. 
The point is after dualizing the dual of $R^n$ maps onto Ext$^1_R(M,R)$.
It turns out one can continue this process and at the next stage obtain
a short exact sequence, $0\ra P\ra N_1\ra M\ra 0$ where now
Ext$^1_R(N_1,R) = \text{Ext}^2_R(N_1,R) = 0$ and $P$ has projective dimension at most
$1$. Continuing until the dimension $d$ of $R$ we eventually obtain a
sequence $0\ra Q\ra C\ra M\ra 0$ where $Q$ has finite projective dimension
and Ext$^1_R(C,R) = \text{Ext}^2_R(C,R) = ... = \text{Ext}^d_R(C,R) = 0$.
If $R$ is assumed to be local and Gorenstein, the duality forces
$C$ to be Cohen-Macaulay, which gives Theorem 11.2. This argument can
be found, with full details, in \cite{EG, (5.5)}.

Gorenstein rings are now part of the basic landscape of mathematics.
 A search for the word `Gorenstein' (algebras or rings) in MathSci Net reveals around
1,000 entries. There are many offshoots which were not mentioned in this article.
The ubiquity paper
and all of its manifestations are indeed ubiquitous! 
\bigskip

\centerline{\bf Bibliography}
\bigskip
\refstyle{A}


\Refs\nofrills{}
\widestnumber\key{AuBu2}

\ref \key{A} \by R. Ap\' ery
\paper La g\' eom\' etrie alg\' ebrique
\jour Bull. Soc. Math. France
\vol 71
\pages 46--66
\yr 1943
\endref

\ref\key{AB}\by M. Auslander and M. Bridger
\paper Stable Module Theory
\jour Mem. Amer. Math. Soc.
\vol 94
\yr 1969
\endref

\ref\key{AuBu1}
\manyby M. Auslander and D. Buchsbaum
\paper Homological dimension in local rings
\jour Trans. Amer. Math. Soc.
\vol 85
\yr 1957
\pages 390--405
\endref

\ref\key{AuBu2}
\bysame
\paper Uniques factorization in regular local rings
\jour Proc. Nat. Acad. Sci. U.S.A.
\vol 45
\yr 1959
\pages 733-734
\endref
\ref\key{ABu}\by M. Auslander and R. Buchweitz
\paper The homological theory of maximal Cohen-Macaulay approximations
\paperinfo preprint
\endref


\ref\key{AvB}\by A. Borek and L. Avramov
\paper Factorial extensions of regular local rings and invariants of finite 
groups
\jour J. reine angew. Math.
\vol 478
\pages 177-188
\yr 1996
\endref

\ref\key{AvF}\by L. Avramov and H.-B. Foxby
\paper Locally Gorenstein homomorphisms
\jour Amer. J. Math.
\vol 114
\pages 1007--1047
\yr 1992
\endref

\ref \key{Ba1} \manyby H. Bass
\paper Injective dimension in Noetherian rings 
\jour Trans. Amer. Math. Soc.
\vol 102
\pages 18--29 
\yr 1962
\endref


\ref\key{Ba2} \bysame
\paper On the ubiquity of Gorenstein rings
\jour Math. Zeitschrift
\vol 82
\pages 8--28 
\yr 1963
\endref

\ref\key{Br}\by H. Bresinsky
\paper Monomial Gorenstein ideals
\jour manu. math.
\vol 29
\pages 159--181
\yr 1979
\endref

\ref\key{BH}
\by W. Bruns and J. Herzog
\book Cohen-Macaulay Rings
\bookinfo Cambridge Studies in Advanced Mathematics
\vol 39
\yr 1993
\publ Cambridge University Press
\endref

\ref\key{BE} \by D. Buchsbaum and D. Eisenbud
\paper Algebra structures for finite free resolutions and some structure
theorems for ideals of codimension $3$
\jour Amer. J. Math.
\vol 99
\pages 447--485
\yr 1977
\endref


\ref\key{Ei}\by D. Eisenbud
\book Commutative Algebra with a View Toward Algebraic Geometry
\bookinfo Graduate Texts in Mathematics
\vol 150
\publ Springer Verlag, New York
\yr 1996
\endref

\ref\key{EGH}
\by D. Eisenbud, M. Green and J. Harris
\paper Cayley-Bacharach theorems and conjectures
\jour Bull. Amer. Math. Soc.
\vol 33
\pages 295--324
\yr 1996
\endref

\ref
\key{EG}
\by G. Evans and P. Griffith
\book Syzygies
\publ Cambridge University Press
\bookinfo London Math. Soc. Lecture Notes
\vol 106
\yr 1985
\endref

\ref\key{Ga}\by E. Garfield
\paper Journal Citation Studies 36. Pure and Applied Mathematics Journals:
What they cite and vice versa
\jour Current Contents
\vol 5(15):5-10
\yr 1982
\pages 34--41
\endref

\ref\key{Go} \by D. Gorenstein
\paper An arithmetic theory of adjoint plane curves
\jour Trans. Amer. Math. Soc.
\vol 72
\pages 414-436
\yr 1952
\endref

\ref\key{G}\by P. Griffith
\paper Normal extensions of regular local rings
\jour J. Algebra
\vol 106
\pages 465--475
\yr 1987
\endref

\ref\key{Gro} \by W. Gr\" obner
\paper \" Uber irreduzible Ideale in kommutativen Ringen
\jour Math. Annalen
\vol 110
\yr 1934
\pages 197--222
\endref

\ref\key{Gr} \by A. Grothendieck
\paper Th\' eor\` emes de dualit\' e pour les faiceaux alg\' ebriques
coh\' erents
\jour S\' eminaire Bourbaki
\vol 9 (149)
\yr 1956/57
\endref

\ref\key{HeK1}
\manyby J. Herzog and E. Kunz
\paper Die Wertehalbgruppe eines lokalen Rings der Dimension $1$
\jour Ber. Heidelberger Akad. Wiss.
\vol 2
\yr 1971
\endref

\ref\key{HeK2}
\bysame
\book Der kanonische Modul eines Cohen-Macaulay Rings
\bookinfo Lectures Notes in Math.
\vol 238
\publ Springer-Verlag
\yr 1971
\endref

\ref\key {HM}\by J. Herzog and M. Miller
\paper Gorenstein ideals of deviation two
\jour  Comm. Algebra
\vol  13
\yr 1985
\pages 1977--1990
\endref



\ref \key{Ho1}
\manyby M. Hochster
\book Topics in the homological theory of modules over commutative rings
\publ C.B.M.S. Regional Conf. Ser. in Math. No. {\bf 24}
\publaddr A.M.S., Providence, R.I.
\yr 1975
\endref

\ref \key {Ho2}
\bysame
\paper Computing the canonical module of a ring of invariants
\jour Cont. Math.
\vol 88
\pages 43--83
\yr 1989
\endref

\ref\key{Ho3}\bysame
\paper Cyclic purity versus purity in excellent Noetherian rings
\jour Trans. Amer. Math. Soc. \vol 231  \yr 1977 \pages 463--488
\endref


\ref
\key {HoR}
\by Hochster, M. and J.L. Roberts
\paper Rings of invariants of reductive groups acting on regular rings are Cohen-Macaulay
\jour Advances in Math.
\vol 13
\yr 1974
\pages 115--175
\endref

\ref
\key{HU1}
\manyby C. Huneke and B. Ulrich
\paper The structure of linkage
\jour Annals Math.
\vol 126
\pages 277--334
\yr 1987
\endref

\ref\key{HU2}\bysame
\paper Divisor class groups and deformations
\jour  American  J.  Math.  
\vol 107 
\yr 1985
\pages 1265--1303
\endref

\ref\key {Ia}
\by A. Iarrobino
\paper Associated graded algebra of a Gorenstein Artin algebra
\jour Memoirs Amer. Math. Soc.
\vol 107
\yr 1994
\endref

\ref\key{Kn}\by F. Knop
\paper Der kanonische Modul eines Invariantenrings
\jour J. Alg.
\vol 127
\pages 40--54
\yr 1989
\endref
\ref
\key {Ku}
\by E. Kunz
\paper Almost complete intersections are not Gorenstein rings
\jour J. Alg.
\vol 28
\yr 1974
\pages 111--115
\endref

\ref\key{KM}\by A. Kustin and M. Miller
\paper Structure theory for a class of grade four Gorenstein ideals
\jour Trans. Amer. Math. Soc.
\vol 270
\yr 1982
\pages 287--307
\endref

\ref\key{Mac}
\by F. S. Macaulay
\paper
On the resolution of a given modular system into primary systems
including some properties of Hilbert functions
\jour Math. Annalen
\yr 1913
\vol 74
\pages 66--121
\endref


\ref\key{Ma} \by E. Matlis
\paper Injective modules over noetherian rings
\jour Pacific J. Math.
\vol 8
\pages 511--528
\yr 1958
\endref

\ref\key{Mu}
\by M.P. Murthy
\paper A note on factorial rings
\jour Arch. Math.
\vol 15
\pages 418--420
\yr 1964
\endref

\ref\key{No}\by E. Noether
\paper Idealtheorie in Ringbereichen
\jour Math. Annalen
\vol 83
\pages 24--66
\yr 1921
\endref


\ref\key{NR} \by D.G. Northcott and D. Rees
\paper Principal systems
\jour Quart. J. Math.
\vol 8
\pages 119--127
\yr 1957
\endref


\ref
\key{PS1} \manyby C. Peskine and L. Szpiro
\paper Dimension projective finie et cohomologie locale
\jour I.H.E.S.
\vol 42
\pages
\yr 1973
\endref

\ref
\key{PS2} \bysame 
\paper Liaison des vari\'et\'es alg\'ebriques
\jour Invent. Math. 
\vol 26 
\yr 1974
\pages 271--302
\endref



\ref
\key {Rob1}
\manyby P. Roberts
\paper Le th\'eor\`eme d'intersection
\jour C. R. Acad. Sc. Paris S\'er. I
\vol 304
\yr 1987
\pages 177--180
\endref

\ref
\key{Rob2}
\bysame
\paper Rings of type $1$ are Gorenstein
\jour Bull. London Math. Soc.
\vol 15
\yr 1983
\pages 48--50
\endref



\ref\key{Ro} \by M. Rosenlicht
\paper Equivalence relations on algebraic curves
\jour Annals Math.
\vol 56
\pages 169--191
\yr 1952
\endref

\ref\key{Sa}
\by P. Samuel
\paper Singularit\' es des vari\' et\' es alg\' ebriques
\jour Bull. Soc. Math. de France
\vol 79
\pages 121--129
\yr 1951
\endref
\ref\key{Se1}\manyby J.P. Serre
\paper Sur les modules projectifs
\jour S\' eminaire Dubreil
\yr 1960-61
\endref

\ref\key{Se2}
\bysame
\paper Sur la dimension homologique des anneaux et des modules noeth\' eriens
\jour Proc. Int. Symp. Tokyo-Nikko
\yr 1956
\publ Science Council of Japan
\pages 175--189
\endref

\ref\key{Se3}
\bysame
\paper Algebraic Groups and Class Fields
\jour Graduate Texts in Mathematics
\vol 117
\yr 1988
\publ Springer-Verlag
\endref

\ref\key{St1}
\manyby R. Stanley
\paper Hilbert functions of graded algebras
\jour Advances Math.
\vol 28
\yr 1978
\pages 57--83
\endref

\ref\key{St2}
\bysame
\paper Invariants of finite groups and their applications to combinatorics
\jour Bull. Amer. Math. Soc.
\vol 1
\yr 1979
\pages 475--511
\endref

\ref\key{St3}
\bysame
\book Combinatorics and Commutative Algebra
\bookinfo Progress in Mathematics
\vol 41
\yr 1983
\publ Birkh\" auser
\endref

\ref\key{U}
\by B. Ulrich
\paper Gorenstein rings as specializations of unique factorization domains
\jour J. Alg.
\vol 86
\pages 129--140
\yr 1984
\endref

\ref\key {VaVi}
\by W.V.  Vasconcelos and R. Villarreal
\paper On Gorenstein ideals of codimension four
\jour Proc. Amer. Math. Soc. 
\vol 98 
\yr 1986
\pages 205--210
\endref


\ref\key{Wa}
\by K. Watanabe
\paper Certain invariant subrings are Gorenstein, I, II
\jour Osaka J. Math.
\vol 11
\pages 1-8, 379--388
\yr 1974
\endref
\endRefs
\enddocument